\documentclass{article}
\usepackage[english]{babel}
\usepackage[letterpaper,top=2cm,bottom=2cm,left=3cm,right=3cm,marginparwidth=1.75cm]{geometry}
\pdfoutput=1
\usepackage{amssymb}
\usepackage{appendix}
\usepackage{latexsym}
\usepackage{mathrsfs}
\usepackage{extarrows}
\usepackage{graphicx}
\usepackage{supertabular}
\usepackage{caption}
\usepackage{subfigure}
\usepackage{setspace}
\usepackage{epstopdf}

\newtheorem{theo}{Theorem}[section]

\newtheorem{conj}{Conjecture}
\newtheorem{lemma}[theo]{Lemma}
\newtheorem{claim}[theo]{Claim}

\newtheorem{fact}[theo]{Fact}
\newtheorem{defi}[theo]{Definition}

\def\q{\hspace*{\fill}$\Box$\medskip}

\def\endproofbox{\hskip 1.3em\hfill\rule{6pt}{6pt}}

{%
\noindent{\it Proof}
}%
{%
 \quad\hfill\endproofbox\vspace*{2ex}
}

\title{Tree Embeddings and Tree-Star Ramsey Numbers}
\author{Zilong Yan \thanks{ School of Mathematics, Hunan University, Changsha 410082, P.R. China. Email: zilongyan@hnu.edu.cn.} \and Yuejian Peng \thanks{ Corresponding author. School of Mathematics, Hunan University, Changsha, 410082, P.R. China. Email: ypeng1@hnu.edu.cn. \ Supported in part by National Natural Science Foundation of China (No. 11931002).}
}

\begin{document}
\title{Tree Embeddings and Tree-Star Ramsey Numbers}
\author{Zilong Yan \thanks{ School of Mathematics, Hunan University, Changsha 410082, P.R. China. Email: zilongyan@hnu.edu.cn.} \and Yuejian Peng \thanks{ Corresponding author. School of Mathematics, Hunan University, Changsha, 410082, P.R. China. Email: ypeng1@hnu.edu.cn. \ Supported in part by National Natural Science Foundation of China (No. 11931002).}
}

\maketitle

\begin{abstract}
 We say that a graph $F$ can be embedded into a graph $G$ if $G$ contains an isomorphic copy of $F$ as a subgraph. Guo and Volkmann \cite{GV}  conjectured that if $G$ is a connected graph with at least $n$ vertices and minimum degree at least $n-3$, then any tree with $n$ vertices and maximum degree at most $n-4$ can be embedded into $G$.
In this paper, we give a  result slightly stronger than this conjecture and  obtain a sufficient and necessary condition that a tree with $n$ vertices and  maximum degree at most $n-3$ can be embedded into a connected  graph G with at least $n$ vertices and minimum degree at least $n-3$. Our result implies that the conjecture of Guo and Volkmann is true with one exception.  We also give an application to the Ramsey number of a tree versus a star.
\end{abstract}

Keywords: Tree, Embedding, Ramsey number, Tree-Star Ramsey number.

\section{Introduction}

All graphs considered throughout the paper are simple graphs, i.e. without loops and multiple edges. Let $V(G)$ denote the vertex set of $G$ and let $E(G)$ denote the edge set of $G$. For $v\in V(G)$, let $N(v)=\{u\in V(G)|uv\in E(G)\}$, $N[v]=N(v)\cup\{v\}$, and $d(v)=|N(v)|$. For $S\subseteq V(G)$, denote $N(S)=\cup_{v\in S}N(v)$ and $N[S]=N(S)\cup S$. Let $\delta(G)=min\{d(v)| v\in V(G)\}$ and $\Delta(G)=max\{d(v)| v\in V(G)\}$. Let $K_n$ denote the complete graph on n vertices and $K_{1, m}$ denote the star with $m+1$ vertices. For graphs $G$ and $H$, the {\em Ramsey Number} $R(G, H)$ is the smallest integer $N$ such that any red-blue-coloring of $E(K_N)$ yields a red $G$ or a blue $H$. For graphs $F$ and $G$, we say that an injection $\phi: V(F)\rightarrow V(G)$ is an {\em embedding} of $F$ into $G$ if for any edge $xy$ in $F$, $\phi(x)\phi(y)$ is an edge in $G$. We say that  $F$ {\it can be embedded into} $G$ if there is an  embedding of $F$ into $G$.

Degree conditions for tree embedding have been  studied actively. A well-known conjecture of Erd\H{o}s-S\'os states that any tree with $n$ vertices can be embedded into a simple graph  with average degree exceeding $n-2$. In \cite{HRSW}, Havet, Reed, Stein and Wood studied the conditions on the maximum degree and the minimum degree of a graph to embed any tree with $n$ vertices, they proposed an interesting conjecture that any tree with $n$ vertices can be embedded into a simple graph with maximum degree more than $n+1$ and minimum degree at least $\lfloor {2(n+1) \over 3}\rfloor$. The Loebl-Koml\'os-S\'os conjecture \cite{EFLS95} states that  any tree with $n$ vertices can be embedded into a simple graph with median degree at least $n+1$. We study minimum degree conditions to embed trees, and it has natural applications on tree-star Ramsey numbers.

In \cite{parson}, Parson determined the path-star Ramsey numbers. A key ingredient in the proof is that $P_n$, a path with $n$ vertices, can be embedded into a graph with minimum degree $n-1$. This is generalized to the following well-known result (see  \cite{cockayne}). % $T_n$  denote a tree on $n$ vertices

\begin{lemma}\label{lemma}
If $G$ is a graph with minimum degree $n-1$, then any tree with $n$ vertices can be embedded into $G$.
\end{lemma}

This lemma implies the Ramsey result of Burr \cite{Burr}: $R(T_n, K_{1, m})\leq m+n-1$, where $T_n$ is a tree with $n$ vertices and $K_{1, m}$ is the star with $m+1$ vertices.
Applying Lemma \ref{lemma}, Cockayne \cite{cockayne} improved the upper bound of $R(T_n, K_{1, m})$ to $m+n-2$ for a special class of trees with $n$ vertices and some values of $m$ and $n$. Further, Guo and Volkmann  \cite{GV} generalized the result of  Cockayne to any tree with $n$ vertices. A key making the generalization possible is that Guo and Volkmann showed that any tree with $n$ vertices other than $K_{1, n-1}$ can be embedded into a connected graph with at least $n$ vertices and minimum degree $n-2$.
They remarked that there are connected graphs $G$ with at least $n$ vertices and minimum degree $n-3$, and trees with $n$ vertices and maximum degree  $n-3$ which cannot be embedded into  $G$. And they proposed the following conjecture.
\begin{conj}\label{gvcon}{\em (\cite{GV})}
If $G$ is a connected graph with at least $n$ vertices and minimum degree at least $n-3$, then any tree  with $n$ vertices and maximum degree at most $n-4$ can be embedded into $G$.
\end{conj}
We show a  result slightly stronger than this conjecture and  obtain a sufficient and necessary condition that trees with $n$ vertices and
maximum degree at most $n-3$ can be embedded into a graph $G$ with minimum degree at least $n-3$. %Applying it, we improve the upper bound of $R(T_n, K_{1, m})$ to $n+m-3$ for some cases of $m$ and $n$. We also obtain the exact value of  $R(T_n, K_{1, m})=n+m-3$ for some cases of $m$ and $n$.

\begin{figure}[ht]
    \centering
    \includegraphics[height=4.5cm]{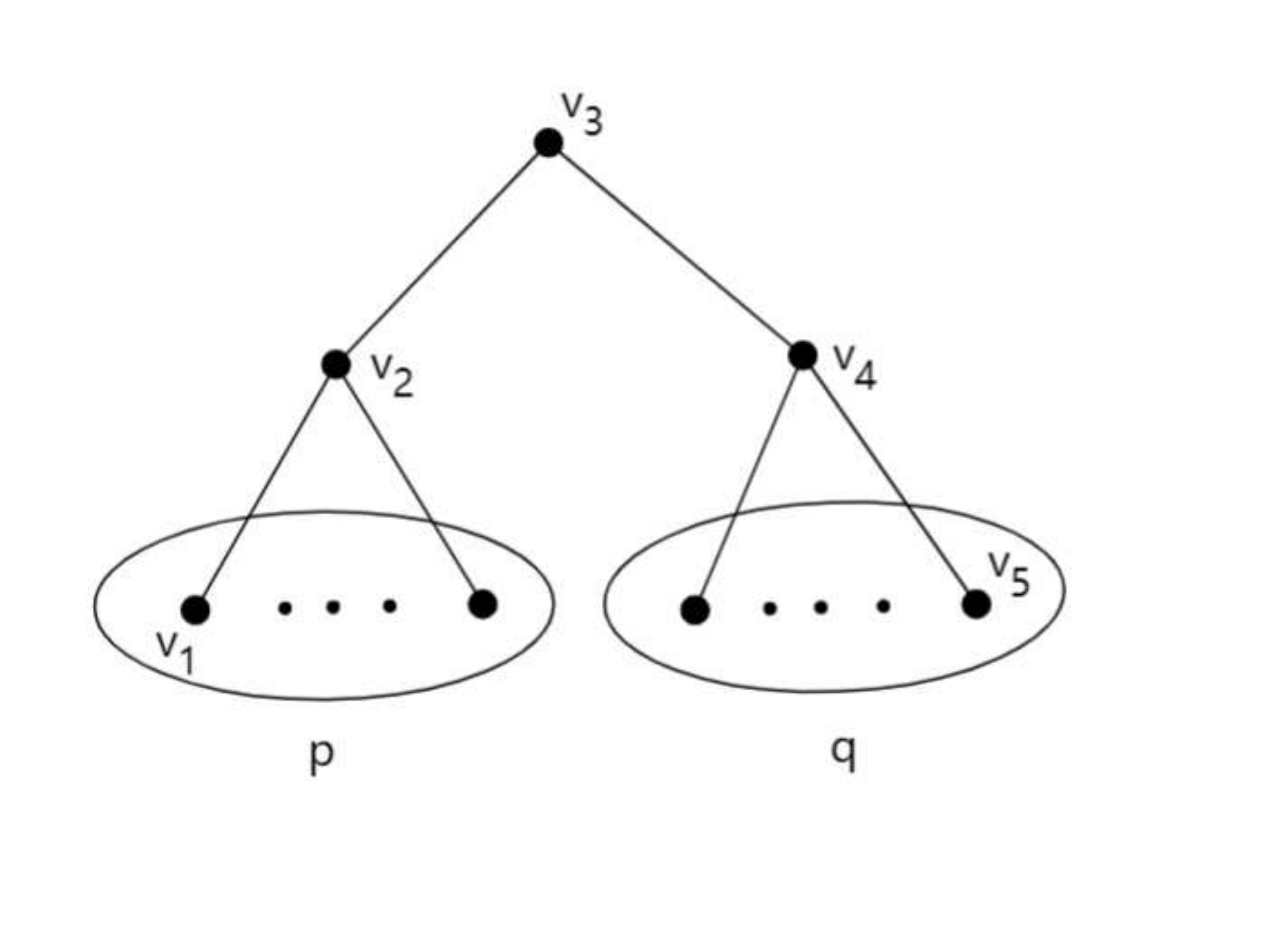}
    \caption{$T(p, q)$}\label{thm}
\end{figure}

Let $p$ and $q$ be positive integers. Let $T(p, q)$ be a tree with a longest path $v_1v_2v_3v_4v_5$ satisfying that $d(v_3)=2$, and $v_2$ and $v_4$ have $p$ and $q$ leaves respectively (see Figure \ref{thm}).
The following is the main result in this paper.
\begin{theo}\label{theo1}
Let $G$ be a connected graph with at least $n$ vertices and minimum degree at least $ n-3$. Let $T_n$ be a tree with $n$ vertices and $\Delta(T_n)\leq n-3$. Then $T_n$ can be embedded into $G$ if and only if the following situations do not happen.\\
(i) $G=K_{n-3, n-3}$ and $T_n=T(p, n-3-p)$, where $p\geq 1$.\\
(ii) $G$ is a balanced complete $(k+1)$-partite graph with $a$ vertices in each part and $ka=n-3$, and $T_n=T(1, n-4)$, where $k\ge 2$ is a positive integer.
\end{theo}
This result implies that Conjecture \ref{gvcon} holds except in situation (i) in Theorem \ref{theo1}.
A non-negative integer $N$ is called a linear combination of two positive integers $p$ and $q$ if there exist non-negative integers $k$ and $l$ such that $N=kp+lq$.
 Applying Theorem \ref{theo1}, we are going to prove the following tree-star Ramsey numbers in Section \ref{sec4}.
\begin{theo}\label{cor1} Let $T_n$ be a tree with $n$ vertices, $\Delta(T_n)\leq n-3$, and $T_n\not= T(p, n-3-p)$. If $m+n-3$ is not a linear combination of $n-1$ and $n-2$, then $R(T_n, K_{1, m})\leq m+n-3$.
\end{theo}

\begin{theo}\label{cor2}
Let $T_n$ be a tree with $n$ vertices and $\Delta(T_n)\leq n-3$. Let $m=k(n-1)+3$ and $0\leq k \leq n-5$.\\% $R(T_n, K_{1, m})\leq m+n-3$ if the following conditions are satisfied \\
(i) If $T_n\not=T(p, n-3-p)$, then $R(T_n, K_{1, m})=m+n-3$.\\
(ii) If $m+n-3$ is not a linear combination of $n-1$, $n-2$, $n-3+a_1$, ..., $n-3+a_d$, where $a_1, a_2, ..., a_d$ are positive divisors of $n-3$ no less than $3$, then $R(T(1, n-4), K_{1, m})=m+n-3$. Otherwise $R(T(1, n-4), K_{1, m})=m+n-2$.\\
(iii) Let $p\geq 2$. If $m+n-3$ is not a linear combination of $2n-6$, $n-1$ and $n-2$, then $R(T(p, n-3-p), K_{1, m})=m+n-3$. Otherwise $R(T(p, n-3-p), K_{1, m})=m+n-2$.
\end{theo}

 The conjecture of Erd\H{o}s-S\'os states that  any tree with $n$ vertices  can be embedded into a graph  with average degree exceeding $n-2$. It is easy to show that a connected graph  with average degree greater than $n-2$ contains a connected subgraph with minimum degree greater than $\frac{1}{2}(n-2)$. So if one can characterize  what kind of trees with $n$ vertices can be embedded into a connected graph $G$ with at least $n$ vertices and $\delta(G)\geq \frac{1}{2}(n-1)$, then these trees will satisfy the conjecture.   In general, it is interesting to study what kind of $T_n$ can be embedded into a connected graph $G$ with at least $n$ vertices and $\delta(G)\geq t$.

In Section \ref{sec2}, we give some crucial lemmas in the proof of Theorem \ref{theo1}. The proof of Theorem \ref{theo1} will be given in Section \ref{sec3}. In Section \ref{sec4}, we prove Theorems \ref{cor1} to \ref{cor2}.

\section{Preparations}\label{sec2}
\begin{defi}
A labelling $\{v_1, v_2, ..., v_n\}$ of vertices of a tree with $n$ vertices is called a conventional labelling if for each $j\in [2, n]$, $|N(v_j)\cap\{v_1, ..., v_{j-1}\}|=1$. To simplify the notation, we always denote the unique vertex in $N(v_j)\cap\{v_1, ..., v_{j-1}\}$ by $v_{j'}$ for a conventional labelling.
\end{defi}

Note that a tree $T$ always has a conventional labelling. We may take any vertex in $V(T)$ as $v_1$ and order the other vertices as $v_2, ..., v_n$ in increasing order of the distances to $v_1$.

For given graphs $F$ and $G$, we say that {\em  a subset $U\subseteq V(F)$ has been embedded into $G$} (or $\phi$ is an {\em  embedding of $U$ to G}) if an embedding $\phi$ of $F[U]$ into $G$ has been established. If $U\subseteq V(F)$ has been embedded to $G$ via an embedding $\phi$, we say that $\phi$ can be extended to a vertex $v\in V(F)\setminus U$ if there exists $v'\in V(G)$ such that the extension function $\hat\phi: U\cup\{v\}\rightarrow V(G)$ defined by $\hat\phi(u)=\phi(u)$ for $u\in U$ and $\hat\phi(v)=v'$ is an embedding of $F[U\cup\{v\}]$ into $G$.

\begin{lemma}\label{lemma23}
Let $G$ be a connected graph with at least $n$ vertices. Let $\{v_1, v_2, ..., v_n\}$ be a conventional labelling of $T_n$. Let $j\in \{2, ..., n\}$. Let $\phi$ be an embedding of $\{v_1, ..., v_{j-1}\}$ into $G$.\\
(i) If
\begin{equation}
|N(\phi(v_{j'}))\setminus\{\phi(v_1), ..., \phi(v_{j-1})\}|\geq 1, \label{1}
\end{equation}
then $\phi$ can be extended to $v_j$.\\
(ii) If $\delta(G)\geq n-1$, then $\phi$ can be extended to $v_j$.\\
(iii) Let $\delta(G)\geq n-2$. If $2\leq j\leq n-1$, then $\phi$ can be extended to $v_j$. If $j=n$ and there exists $k\in \{1, 2, ..., n-1\}\setminus\{n'\}$ such that $\phi(v_k)\notin N(\phi(v_{n'}))$, then $\phi$ can be extended to $v_j$.\\
(iv) Let $\delta(G)\geq n-3$. If $2\leq j\leq n-2$, then $\phi$ can be extended to $v_j$. If $j=n-1$ and $|\{\phi(v_1), ..., \phi(v_{n-2})\}\setminus N(\phi(v_{(n-1)'}))|\geq 2$, then $\phi$ can be extended to $v_j$. If $j=n$ and $|\{\phi(v_1), ..., \phi(v_{n-1})\}\setminus N(\phi(v_{n'}))|\geq 3$, then $\phi$ can be extended to $v_j$.
\end{lemma}
{\em Proof of Lemma \ref{lemma23}.} The condition in (i) guarantees that there is a free vertex in $N(\phi(v_{j'}))$ to embed $v_j$, therefore $v_j$ can be embedded to a vertex in $N(\phi(v_{j'}))$, and this proves (i).  (ii), (iii) and (iv) are direct implications of (i).
\q

\begin{lemma}\label{lemma3}(\cite{Dirac})
Let $G$ be a connected graph with minimum degree $\delta(G)$. Then $G$ contains a path of length $min\{n, 2\delta(G)+1\}$.
\end{lemma}

\begin{lemma}\label{lemma25}
Let $G$ be a connected graph with minimum degree $\delta$. Let $S_0\subseteq V(G)$ and let $S_1$ be a proper subset of $S_0$ such that $N(S_1)\subseteq S_0$. If there exists $\alpha\in S_0$ such that $N(\alpha)\not\subset S_0$, then there exists a path of length at least $|S_1|+\delta-|S_0|+1$ starting from a neighbor $\alpha$ outside $S_0$.
\end{lemma}
{\em Proof of Lemma \ref{lemma25}.} Let $\alpha_1\in N(\alpha)\setminus S_0$. Let $\alpha_1\alpha_2...\alpha_g$ be a maximum path outside $S_0$. Since $N(S_1)\subseteq S_0$, $N(\alpha_g)\subseteq S_0\cup\{\alpha_1, \alpha_2, ..., \alpha_{g-1}\}\setminus S_1$. Since $|N(a_g)|\geq \delta$, $|S_0|+g-1-|S_1|\geq \delta$. So $g\geq |S_1|+\delta-|S_0|+1$.
\q

\begin{lemma}\label{lemma26}
Let $G$ be a non-complete and connected graph, then there exist $\alpha, \beta, \gamma\in V(G)$ such that $\alpha\beta, \alpha\gamma\in E(G)$ and $\beta\gamma\notin E(G)$.
\end{lemma}
{\em Proof of Lemma \ref{lemma26}.} Take a shortest path between two non-adjacent vertices, then three consecutive vertices in this path satisfy the condition.
\q

\section{Embedding trees}\label{sec3}
We give the proof of Theorem \ref{theo1} below.

{\em Proof of Theorem \ref{theo1}.} Let $T_n$ be a tree with $n$ vertices. Let $P$ be a longest path of $T_n$. Let $|P|$ be the number of vertices of path $P$ and $i=|P|-2$. Since $|P|=3$ implies that $T_n=K_{1,n-1}$, it contradicts $\Delta(T_n)\leq n-3$. So $|P|\geq 4.$ For $|P|\ge 5$, we label the two end vertices as $v_{n-1}$ and $v_n$, all other vertices in $P$ as $v_1, v_2, ..., v_i$ such that $P=v_{n-1}v_1v_2...v_iv_n$, i.e. $v_{n-1}$ is connected to $v_1$, $v_1$ is connected to $v_2$,  $\cdots$, $v_{i-1}$ is connected to $v_i$, and $v_i$ is connected to $v_n$. We order the vertices outside $P$ as $\{v_{i+1}, ..., v_{n-2}\}$ in non-decreasing order of the distances to $P$. For example, $v_{i+1}$ is connected to a vertex in $P$. Note that $\{v_1, v_2, ..., v_{n-1}, v_n\}$ is a conventional labelling. We will show that there is no embedding $\phi: V(T_n)\rightarrow V(G)$ if and only if situation (i) or (ii) in Theorem \ref{theo1} happens.

By Lemma \ref{lemma}, we can easily embed $\{v_1, ..., v_{n-2}\}$ into $G$. So what we need to work on is to embed $v_1$ and $v_i$ properly  so that there is a free vertex for us to embed $v_{n-1}$ in the neighborhood of the vertex  $v_1$ embedded to, and there is a free vertex for us to embed $v_{n}$ in the neighborhood of the vertex  $v_i$ embedded to.

By Lemma \ref{lemma26}, there exist $u_1$, $u_2$ and $u_3$ in $G$ satisfying that $u_2\in N(u_1)$ and $u_3\in N(u_2)\setminus N(u_1)$. Let $\phi(v_j)=u_j$ for $j\in [3]$. By Lemma \ref{lemma23} (iv), $\phi$ can be extended to $\{v_4, ..., v_{i-1}\}$. Denote $\phi(v_j)=u_j$ for $j\in [i-1]$. Denote $K_i=\{\phi(v_1), ..., \phi(v_{i-1})\}$. %Let $P=v_{n-1}v_1v_2...v_iv_n$.

{\em Case 1.} $|P|\geq 6$.

{\em Case 1.1.}  There exists a vertex, denoted by $u_i$, in $N(u_{i-1})\setminus K_i$ such that $|K_i\setminus N(u_i)|\ge 2$.

In this case, let $\phi(v_i)=u_i$. Since $\delta(G)\geq n-3$, $u_1, u_3\notin N(u_1)$, by Lemma \ref{lemma23} (iv), $\phi$ can be extended to $\{v_1, ..., v_{n-1}\}$. Further  $|K_i\setminus N(u_i)|\ge 2$ implies that $|\{\phi(v_1), ..., \phi(v_{n-1})\}\setminus N(u_i)|\geq 3$, by Lemma \ref{lemma23} (iv), there exists an embedding from $T_n$ into $G$.

{\em Case 1.2.}  There exists a vertex, denoted by $u_i$, in $N(u_{i-1})\setminus K_i$ such that $|K_i\setminus N(u_i)|=1$, and $|K_i\setminus N(u)|\le1$ for any $u\in N(u_{i-1})\setminus K_i$.

Note that $|K_i\setminus N(u_i)|=1$ implies that there exists a vertex $u_x$ in $K_i$ such that $u_xu_i\notin E(G)$.

If $d(u_i)\geq n-2$, let $\phi(v_i)=u_i$. Since $u_3\notin N(u_1)$, by Lemma \ref{lemma23} (iv), $\phi$ can be extended to $\{v_{i+1}, ..., v_{n-1}\}$. Since $u_x, u_i\notin N(u_i)$, $|N(u_i)\setminus\{\phi(v_1), ..., \phi(v_{n-1})\}|\geq 1$ and there is at least one free vertex in $N(u_i)$ to embed $v_n$. So it is sufficient to consider $d(u_i)=n-3$. Let $K=K_i\cup N[u_i]$, then $|K|=n-1$.

\begin{claim}\label{claim0}
$N(u_1)\subseteq K$ if $T_n$ cannot be embedded into G.
\end{claim}
{\em Proof of Claim \ref{claim0}.} Suppose that there exists a vertex denoted by  $u_{n-1}$ in $N(u_1)\setminus K$. Since $u_3\notin N(u_1)$, by Lemma \ref{lemma23} (iv), $\phi$ can be extended to $\{v_1, ..., v_{n-1}\}$. Furthermore, we can guarantee that some vertex $v_a\in T_n\setminus\{v_n\}$ is embedded to $u_{n-1}$, where $a\in [i+1, n-1]$ (This is always possible since we can always embed $v_{n-1}$ to $u_{n-1}$ if $u_{n-1}$ is not previously used). Now $u_x, u_i, u_{n-1}\notin N(u_i)$, then there exists a free vertex in $N(u_i)\setminus\{\phi(v_1), ..., \phi(v_{n-2}), \phi(v_{n-1})\}$ to embed $v_n$ and $T_n$ can be embedded into $G$.
\q

\begin{claim}\label{claim1}
$N(u_2)\subseteq K$ if $T_n$ cannot be embedded into G.
\end{claim}
{\em Proof of Claim \ref{claim1}.} If there exists a vertex $u_1'\in N(u_2)\setminus K$, then reassign $\phi(v_1)=u_1'$ and let $\phi(v_i)=u_i$. By the choice of $u_1'$, we know that $u_i\notin N(\phi(v_1))$. By Lemma \ref{lemma23} (iv), $\phi$ can be extended to $\{v_1, ..., v_{n-1}\}$. Recall that there exists a vertex $u_x\in K_i$ such that $u_xu_i\notin E(G)$. We claim that $u_x\not=u_1$. Since $d(u_1)\geq n-3$, $u_1, u_3\notin N(u_1)$, $N(u_1)\subseteq K$ and $|K|=n-1$, $u_1u_i\in E(G)$. Since $u_x, u_1'\notin N(u_i)$, there is at least one free vertex in $N(u_i)\setminus\{\phi(v_1), ..., \phi(v_{n-1})\}$ to embed $v_n$. A contradiction.
\q

\begin{claim}\label{claim2}
$N(u_t)\subseteq K$ for $1\leq t\leq i-1$ if $T_n$ cannot be embedded into G.
\end{claim}
{\em Proof of Claim \ref{claim2}.} Use induction on t. For $t=1$ or 2, it is guaranteed by Claim \ref{claim1}. Let $3\leq t\leq i-1.$ We assume that $N(u_s)\subseteq K$ for $1\leq s\leq t-1$. If $N(u_t)\not\subset K$, then there exists a vertex $u_{t-1}'\in N(u_t)\setminus K$. By Lemma \ref{lemma25}, there exists a path $u_{t-1}'u_{t-2}'...u_1'$ outside K. Reassign $\phi(v_{\alpha})=u_{\alpha}'$ for $1\leq\alpha\leq t-1$ and let $\phi(v_i)=u_i$. Since $u_{\alpha}'\notin N(u_i)$ for $1\leq\alpha\leq t-1$ and $t\geq 3$, by Lemma \ref{lemma23} (iv), $\phi$ can be extended to $\{v_{i+1}, ..., v_n\}$.
\q

\begin{claim}\label{claim3}
$N(K)\subseteq K$ if $T_n$ cannot be embedded into G.
\end{claim}
{\em Proof of Claim \ref{claim3}.} If there exists a vertex $u_{i-1}'\in N(u_i)\setminus K_i$ such that there exists a vertex $u_{i-2}'\in N(u_{i-1}')\setminus K$, by Lemma \ref{lemma25}, there exists a path $u_{i-2}'u_{i-3}'...u_1'$ outside K. Reassign $\phi(v_{\alpha})=u_{\alpha}'$ for $1\leq\alpha\leq i-1$ and let $\phi(v_i)=u_i$. Since $i\geq 4$ guarantees that $u_1', u_2'\notin N(u_i)$, by Lemma \ref{lemma23} (iv), $\phi$ can be extended to $\{v_{i+1}, ..., v_n\}$.
\q

By Claim \ref{claim3}, if $T_n$ cannot be embedded into $G$, then $K$ is a component with $n-1$ vertices, contradicting that $G$ is a connected graph with at least $n$ vertices. So $T_n$ can be embedded into $G$ in this case.

{\em Case 1.3.}  For each vertex $u$ in $N(u_{i-1})\setminus K_i$, $|K_i\setminus N(u)|=0$.

Take a vertex, denote by $u_i$, in $N(u_{i-1})\setminus K_i$, and let $\phi(v_i)=u_i$.

If $d(u_i)\geq n-1$, since $u_1, u_3\notin N(u_1)$, by Lemma \ref{lemma23} (iv), $\phi$ can be extended to an embedding on $\{v_1, ..., v_{n-1}\}$. Since $d(u_i)\geq n-1$ and $\phi(v_i)=u_i\notin N(u_i)$, there is at least one free vertex in $N(u_i)\setminus\{\phi(v_1), ..., \phi(v_{n-1})\}$ to embed $v_n$.

If $d(u_i)=n-2$, then $|N[u_i]|=n-1$. We claim that $\cup_{u\in N[u_i]}N[u]=N[u_i]$ if $T_n$ cannot be embedded into $G$. If there exists a vertex, denote by $u_{n-1}$, in $N(u_1)\setminus N[u_i]$, let $\phi(v_i)=u_i$ and extend $\phi$ to an embedding on $\{v_1, ..., v_{n-1}\}$ such that we embed some $v_a\in \{v_{i+1}, ..., v_{n-1}\}$ to $u_{n-1}$ (This is always possible since we can always embed $v_{n-1}$ to $u_{n-1}$ if $u_{n-1}$ is not previously used). Now $u_i, u_{n-1}\notin N(u_i)$, therefore there exists one free vertex in $N(u_i)$ to embedded $v_n$ and $T_n$ can be embedded into $G$. So we have shown that $N[u_1]\subseteq N[u_i]$. Assume that $\cup_{i=1}^t N[u_i]\subseteq N[u_i]$ for $1\leq t\leq i-1$, we show that $N[u_{t+1}]\subseteq N[u_i]$.
If there exists a vertex $u_t'\in N[u_{t+1}]\setminus N[u_i]$, by Lemma \ref{lemma25}, there exists a path $u_t'u_{t-1}'...u_1'$ outside $N[u_i]$. Reassign $\phi(v_\alpha)=u_{\alpha}'$ for $1\leq\alpha\leq t$. Since $u_i\notin N(\phi(v_1))$, by Lemma \ref{lemma23} (iv), $\phi$ can be extended to an embedding on $\{v_1, ..., v_{n-1}\}$. Since $d(u_i)=n-2$ and $\phi(v_i)=u_i$, $\phi(v_a)\notin N(u_i)$ for $1\leq a\leq t$, there is at least one free vertex in $N(u_i)\setminus\{\phi(v_1), ..., \phi(v_{n-1})\}$ to embed $v_n$. Therefore we have shown that $N[K_i]\subseteq N[u_i]$.
If there exist vertices $u', u''$ such that $u'\in N(u_i)\setminus K_i$ and $u''\in N[u']\setminus N[u_i]$, since $d(u_3)\geq n-3$, $|N[u_i]|=n-1$ and $u_1, u_3\notin N(u_3)$ and $N(u_3)\subseteq N[u_i]$, $u_3u'\in E(G)$. Reassign $\phi(v_2)=u'$ and $\phi(v_1)=u''$. Since $u_i\notin N(\phi(v_1))$, by Lemma \ref{lemma23} (iv), $\phi$ can be extended to $\{v_1, v_2, ..., v_{n-1}\}$. Since $u''u_i\notin E(G)$ and $d(u_i)=n-2$, there exists a free vertex in $N(u_i)\setminus\{\phi(v_1), ..., \phi(v_{n-1})\}$ to embed $v_n$ and $T_n$ can be embedded into $G$. So we have shown that $\cup_{u\in N[u_i]}N[u]=N[u_i]$ if $T_n$ cannot be embedded into $G$. Therefore $N[u_i]$ is a component of $n-1$ vertices, contradicting that $G$ is a connected graph with at least $n$ vertices.

Now we discuss $d(u_i)=n-3$. Since $u_1u_3\notin E(G)$ and $d(u_3)\geq n-3$, there exists a vertex, denoted by $u_2'$, in $N(u_3)\setminus N[u_i]$. If there exists a vertex $u_1'\in N(u_2')\setminus N[u_i]$, then reassign $\phi(v_2)=u_2'$ and $\phi(v_1)=u_1'$. Since $\phi(v_1), \phi(v_2)\notin N(\phi(v_i))$, by Lemma \ref{lemma23} (iv), $\phi$ can be extended to $\{v_1, ..., v_n\}$ and $T_n$ can be embedded into $G$. So $N(u_2')\subseteq N[u_i]$. Since $u_2'u_i\notin E(G)$ and $d(u_2')\geq n-3$, $N(u_2')=N(u_i)$. Since $u_2'u_1, u_2'u_3\in E(G)$, we can reassign $\phi(v_2)=u_2'$. Now $u_i$ has a non-neighbor in $K_i$. Note that $u_3\notin N(\phi(v_1))$ still holds. This situation is exactly {\em Case 1.2} we have proved.

{\em Case 2.} $|P|=5$.

\begin{figure}[ht]
\centering
\begin{minipage}[b]{0.45\linewidth}
\includegraphics[height=4.5cm]{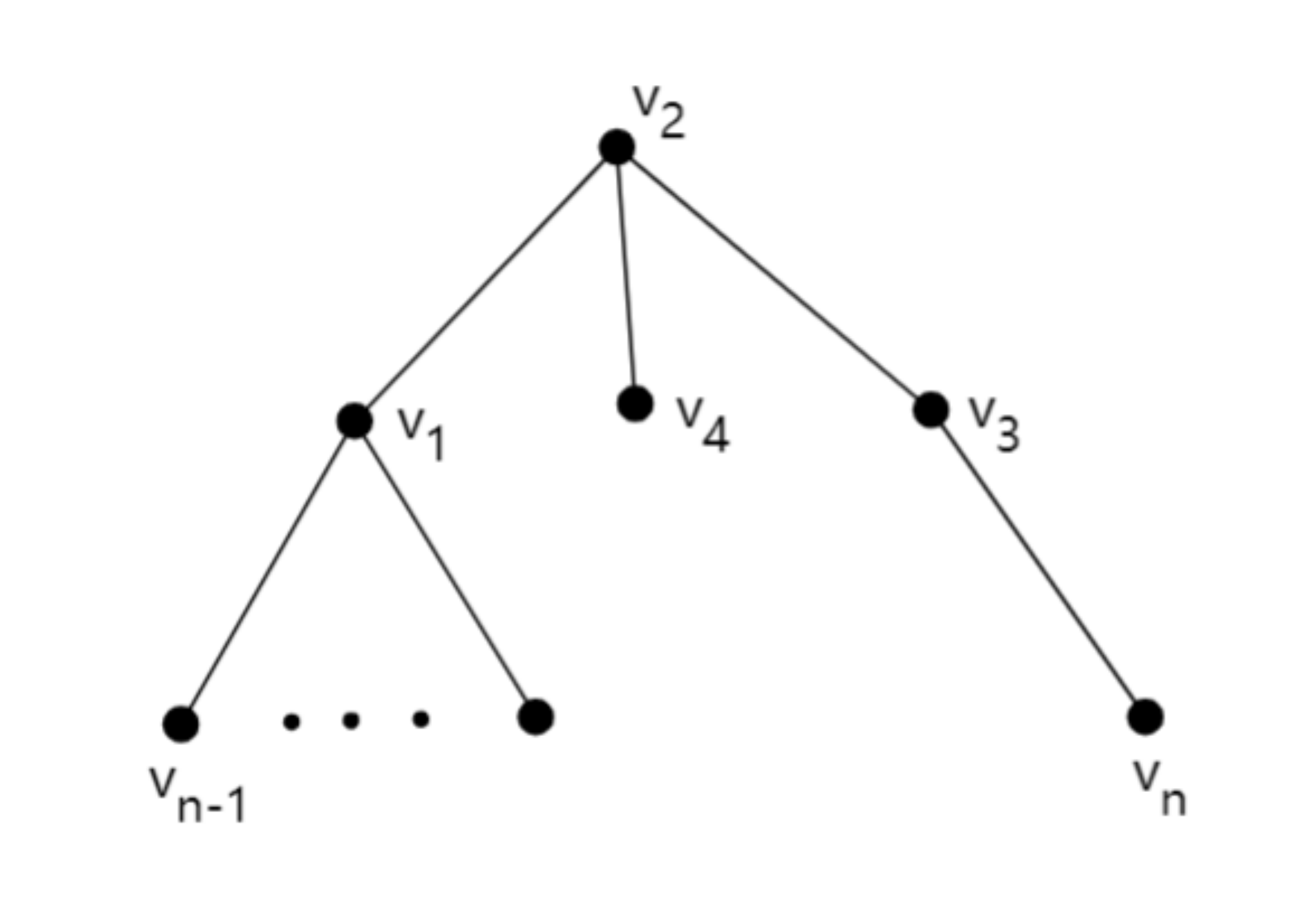}
\caption{Case $d_{T_n}(v_2)\ge 3$}
\label{fig2}
\end{minipage}
\quad
\begin{minipage}[b]{0.45\linewidth}
\includegraphics[height=4.5cm]{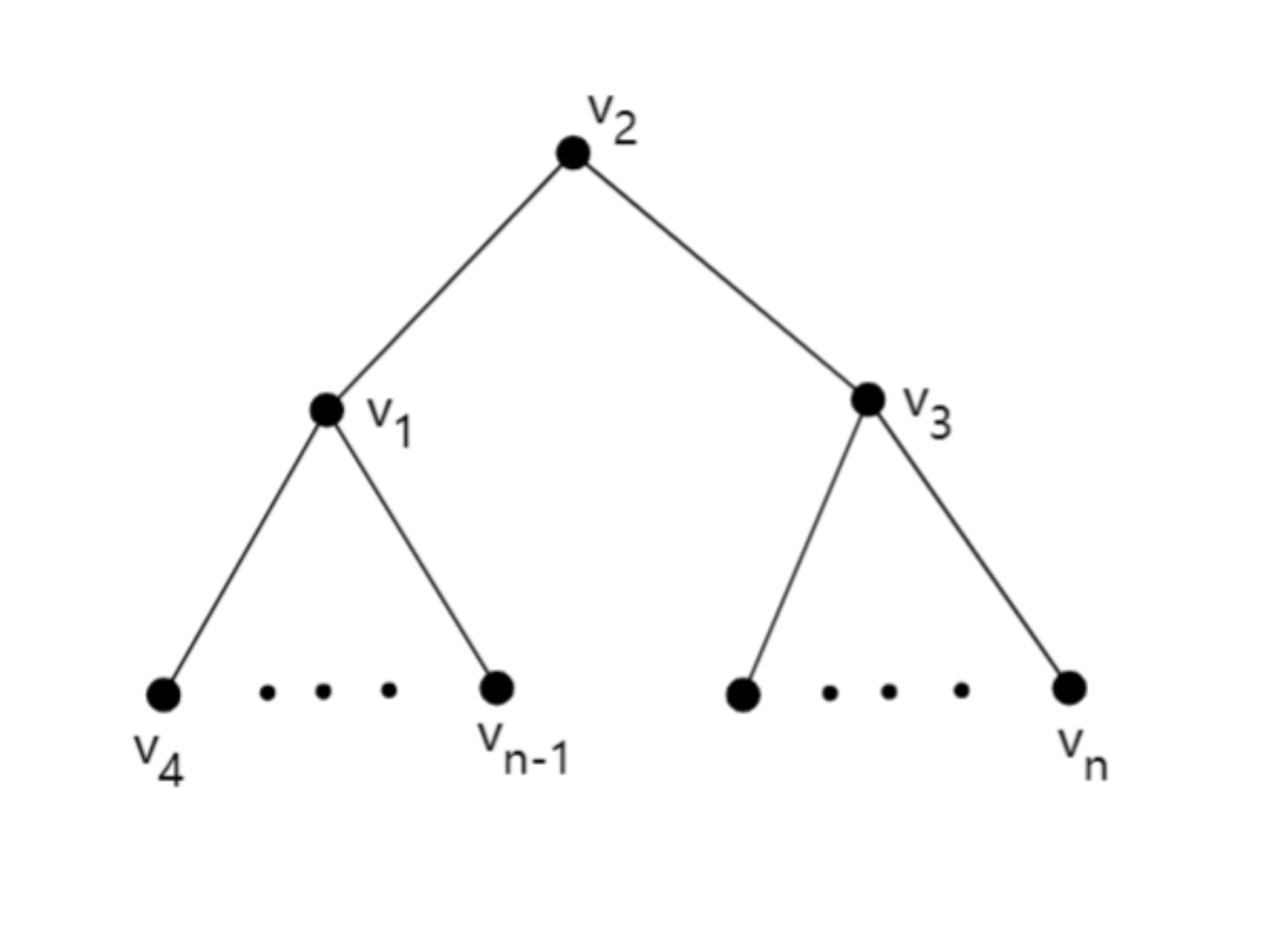}
\caption{Case $d_{T_n}(v_2)=2$}
\label{fig3}
\end{minipage}
\end{figure}
In this case, $P=v_{n-1}v_1v_2v_3v_n$. If $d_{T_n}(v_2)\geq 3$, then there exists a vertex $v_4$ connected to $v_2$ (see Figure \ref{fig2}). If $d_{T_n}(v_2)=2$, then $N_{T_n}(v_2)=\{v_1, v_3\}$. If $T_n=P$, since $|V(G)|\geq5$ and $\delta(G)\geq 2$, by Lemma \ref{lemma3}, $P=T_n$ can be embedded in $G$. So if $d_{T_n}(v_2)=2$, we may assume that $d_{T_n}(v_1)\geq 3$ or $d_{T_n}(v_3)\geq 3$ (see Figure \ref{fig3}).

Let $u_3'$ be a vertex with maximum degree in $V(G)$. If $d(u_3')\geq n-2$, if there exists a vertex $u_2'\in N(u_3')$ such that there exists a vertex $u_1'\in N(u_2')\setminus N[u_3']$, then let $\phi(v_1)=u_1'$, $\phi(v_2)=u_2'$ and $\phi(v_3)=u_3'$. Since $u_1', u_3'\notin N(u_1')$ and $\delta(G)\geq n-3$ and $d(u_3')\geq n-2$, by Lemma \ref{lemma23} (iv), $\phi$ can be extended to $\{v_4, ..., v_n\}$. If there is no such vertex $u_2'$ in $N(u_3')$, then $N[u_3']$ is a component, so $d(u_3')\geq n-1$. If there exists a vertex $u_1'\in N(u_3')$ such that $d(u_1')\geq n-2$, let $\phi(v_1)=u_1'$ and $\phi(v_3)=u_3'$ and $\phi(v_2)\in N(u_1')$, then by Lemma \ref{lemma23} (iv), $\phi$ can be extended to $\{v_4, ..., v_n\}$. So we may assume that $d(u_x)=n-3$ for all $u_x\in N(u_3')$. Since $|N[u_3']|\geq n$, for any vertex $u_3''\in N(u_3')$, we can find vertices $u_1'', u_4''\in N(u_3')\setminus N[u_3'']$. Let $\phi(v_1)=u_1''$ and $\phi(v_2)=u_3'$ and $\phi(v_3)=u_3''$ and $\phi(v_4)=u_4''$. Since $\phi(v_1), \phi(v_3)\notin N(\phi(v_1))$, $d(\phi(v_1))\geq n-3$ and $\phi(v_1), \phi(v_3), \phi(v_4)\notin N(\phi(v_3))$, by Lemma \ref{lemma23} (iv), $\phi$ can be extended to $\{v_4, ..., v_n\}$ and $T_n$ can be embedded into $G$. Therefore, it is sufficient to consider that $G$ is an ($n-3$)-regular connected graph.

{\em Case 2.1.} $G$ is an ($n-3$)-regular connected graph with at least $n$ vertices and $d_{T_n}(v_2)\geq 3$ (see Figure \ref{fig2}).

\begin{figure}[ht]
\centering
\begin{minipage}[b]{0.45\linewidth}
\includegraphics[height=4.5cm]{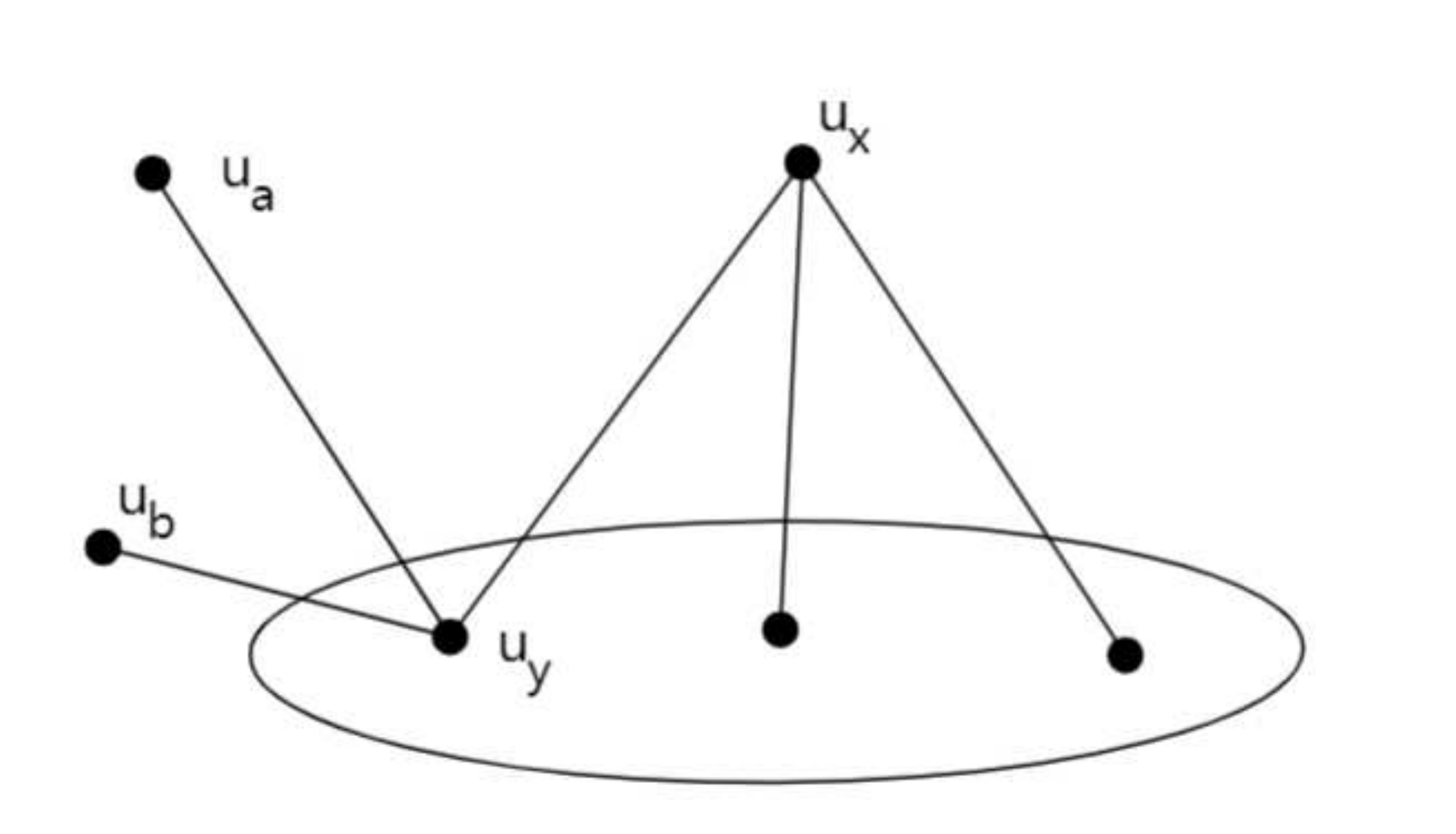}
\caption{Figure for Claim \ref{f2}}
\label{f21}
\end{minipage}
\quad
\begin{minipage}[b]{0.45\linewidth}
\includegraphics[height=4.5cm]{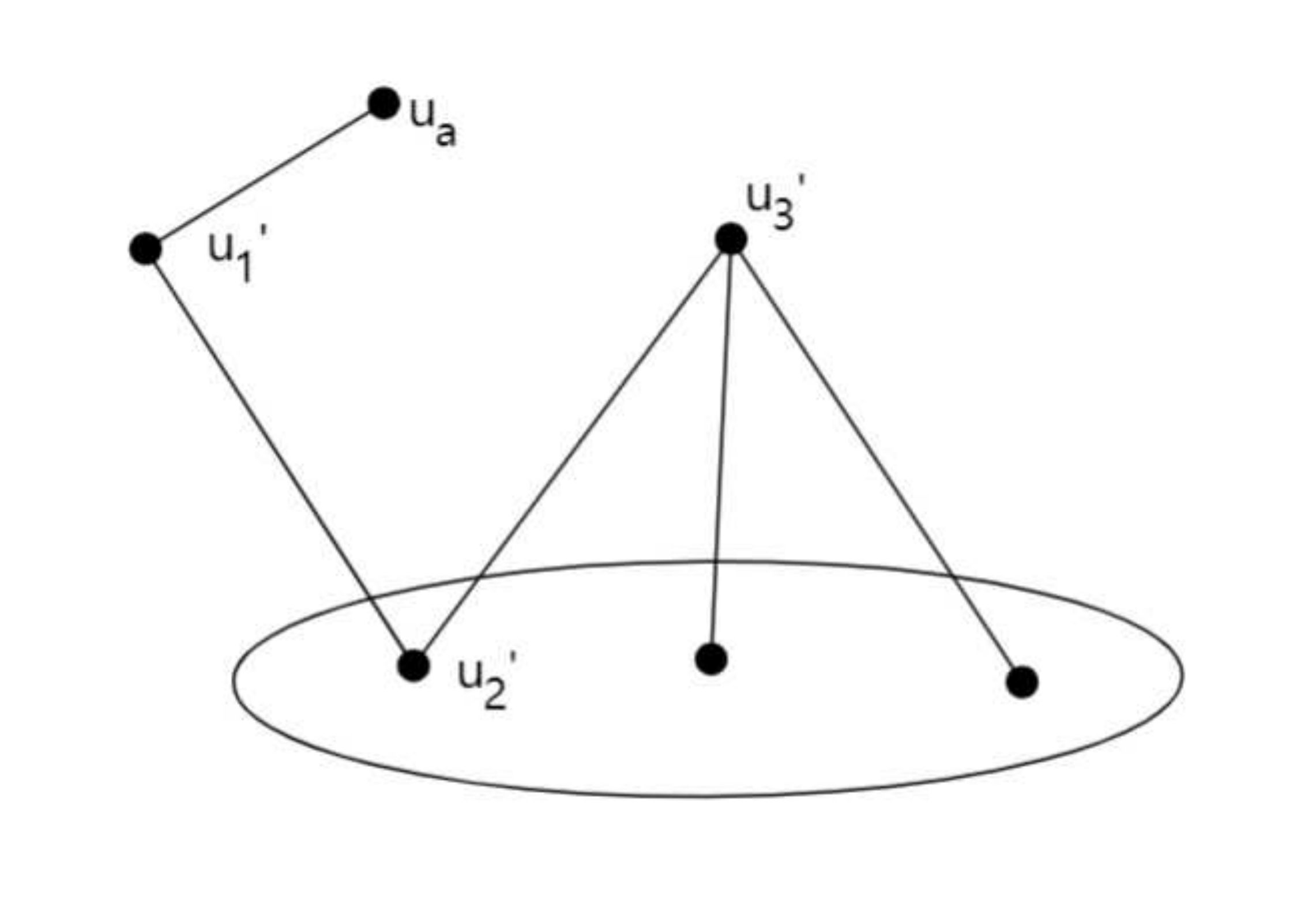}
\caption{Figure for Case 2.1}
\label{f1}
\end{minipage}
\end{figure}

\begin{claim}\label{f2}
Let $G$ be  ($n-3$)-regular and $d_{T_n}(v_2)\geq 3$ (see Figure \ref{fig2}). If there exist vertices $u_x, u_y, u_a, u_b\in V(G)$ such that $u_xu_y\in E(G)$ and $u_a, u_b\in N(u_y)\setminus N[u_x]$ (see Figure \ref{f21}), then $T_n$ can be embedded into G.
\end{claim}
{\em Proof.} Let $\phi(v_1)=u_a$, $\phi(v_2)=u_y$, $\phi(v_3)=u_x$ and $\phi(v_4)=u_b$. Since $G$ is ($n-3$)-regular, $\phi(v_1), \phi(v_2)\notin N(\phi(v_1))$, and $\phi(v_1), \phi(v_4), \phi(v_3)\notin N(\phi(v_3))$, by Lemma \ref{lemma23} (iv), $\phi$ can be extended to $\{v_5, ..., v_n\}$.
\q

Take a vertex $u_3'\in V(G)$. Since $|N[u_3']|=n-2$, $N[u_3']$ is not a component. So there exists vertices $u_1', u_2'$ such that $u_2'\in N(u_3')$ and $u_1'\in N(u_2')\setminus N[u_3']$ (see Figure \ref{f1}). By Claim \ref{f2}, $N[u_2']\subseteq N[u_3']\cup\{u_1'\}$. If $u_1'$ has two neighbors outside $N[u_3']$, then these two neighbors are outside $N[u_2']$. By Claim \ref{f2}, $T_n$ can be embedded into $G$. If $N(u_1')\subseteq N[u_3']$, then $N(u_1')=N(u_3')$. If $T_n$ cannot be embedded into $G$, Claim \ref{f2} implies that $\cup_{u\in N[u_3']}N(u)\cup N(u_1')=N[u_3']\cup\{u_1'\}$. Otherwise there exists a vertex $u\in N(u_3')$ such that $u$ has a neighbor $w$ outside $N[u_3']\cup\{u_1'\}$, then $u$ has two neighbors $w$ and $u_1'$ outside $N[u_3']$, contradicting Claim \ref{f2}. Thus $N[u_3']\cup\{u_1'\}$ is a component on $n-1$ vertices, contradicting that $G$ is a connected graph with at least $n$ vertices. Therefore, there exists exactly one vertex $u_a\in N(u_1')\setminus N[u_3']$. By Claim \ref{f2}, $N(u_1')$ and $N(u_a)$ differ by at most one vertex. So $N(u_1')\subseteq N(u_3')\cup\{u_a\}$. If there exists a vertex $u\in N(u_3')\cap N(u_1')$ connected to $u_a$, then both $u_a$ and $u_1'$ are in $N(u)\setminus N(u_3')$, contradicting Claim \ref{f2}. So $N(u_1')$ and $N(u_a)$ have no common neighbor in $N[u_3']$, contradicting Claim \ref{f2} again. So $T_n$ can be embedded into $G$.

{\em Case 2.2.} $G$ is  ($n-3$)-regular and $d_{T_n}(v_2)=2$ (see Figure \ref{fig3}), i.e. $T_n=T(p, n-3-p)$ for some $p$.

\begin{claim}\label{f3}
Let G be an (n-3)-regular connected graph on at least $n$ vertices. If there exists $u_1', u_2', u_3'\in V(G)$ such that $u_1'u_2', u_2'u_3'\in E(G)$ and $u_1'u_3'\notin E(G)$. Then $T_n$ can be embedded into G or $N(u_1')=N(u_3')$.
\end{claim}
{\em Proof.} Suppose that $N(u_1')\not=N(u_3')$. Since $G$ is regular, there exists a vertex $u_4\in N(u_1')\setminus N[u_3']$. Let $\phi(v_1)=u_1'$, $\phi(v_2)=u_2'$, $\phi(v_3)=u_3'$ and $\phi(v_4)=u_4$. Since $u_1'u_3'\notin E(G)$, and $u_1', u_3', u_4\notin N(u_3')$, by Lemma \ref{lemma23} (iv), $\phi$ can be extended to an embedding of $T_n$.
\q

\begin{figure}[ht]
    \centering
    \includegraphics[height=4.5cm]{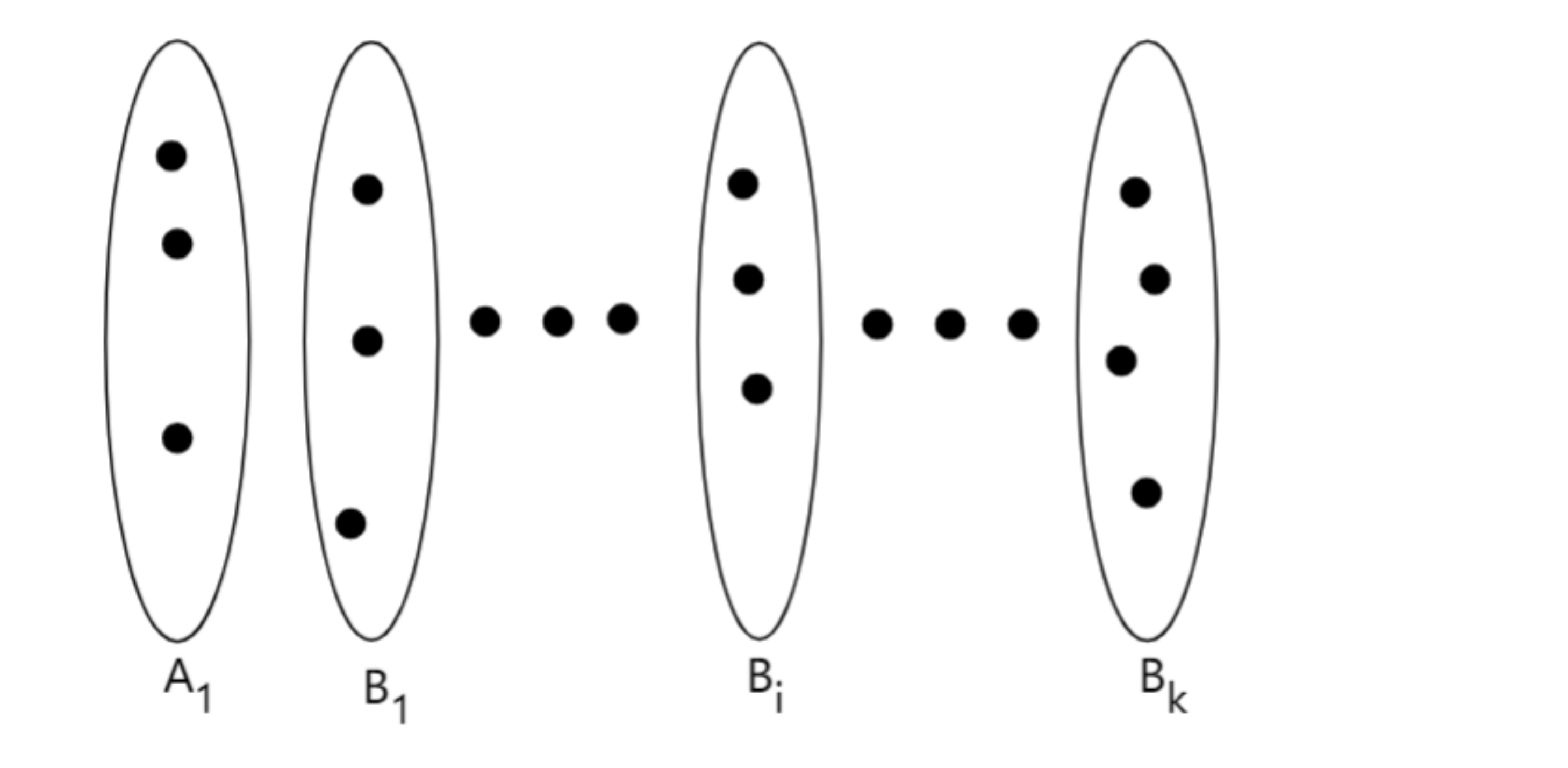}
    \caption{$K_{a, a, \dots, a}$}\label{case2.3}
\end{figure}

Assume that $T_n$ cannot be embedded into $G$. Let $A_1$ be a maximal set in $G$ satisfying that it is an independent set and its any two vertices have a common neighbor. By Claim \ref{f3}, all vertices in $A_1$ have the same set of neighbors, say $B$. Then $|B|=n-3$. We claim that $V(G)=A_1\cup B$. Otherwise, there exists a vertex $u_x\in V(G)\setminus (A_1\cup B)$ such that $u_x$ is connected to a vertex in $B$ and none of the vertices in $A_1$ is adjacent to $u_x$. By Claim \ref{f3}, $N(u_x)=B$. This contradicts the maximality of $A_1$. Since it holds for all such maximal sets, we can divide $B$ into $k$ sets $B_1, B_2, ..., B_k$ such that $B_j$ is such a maximal set for each $j\in [k]$. Since $G$ is regular, $|A_1|=|B_1|=...=|B_k|=a$, $n-3=ak$ and $|V(G)|=n-3+a$ (See Figure \ref{case2.3}). Since $|V(G)|\geq n$, $a\geq 3$.

If $k=1$, since $G$ is ($n-3$)-regular, $G=K_{n-3, n-3}$, and in this case, $T(p, n-3-p)$ cannot be embedded into $K_{n-3, n-3}$. This is because that $T(p, n-3-p)$ is a bipartite graph with 2 and $n-2$ vertices in each of the two parts. This is situation (i) in Theorem \ref{theo1}.

If $k\geq 2$, then let $\phi(v_1)\in B_1$, $\phi(v_2)\in B_2$ and $\phi(v_3)\in A_1$.
\begin{figure}[ht]
\centering
\begin{minipage}[b]{0.45\linewidth}
\includegraphics[height=4.5cm]{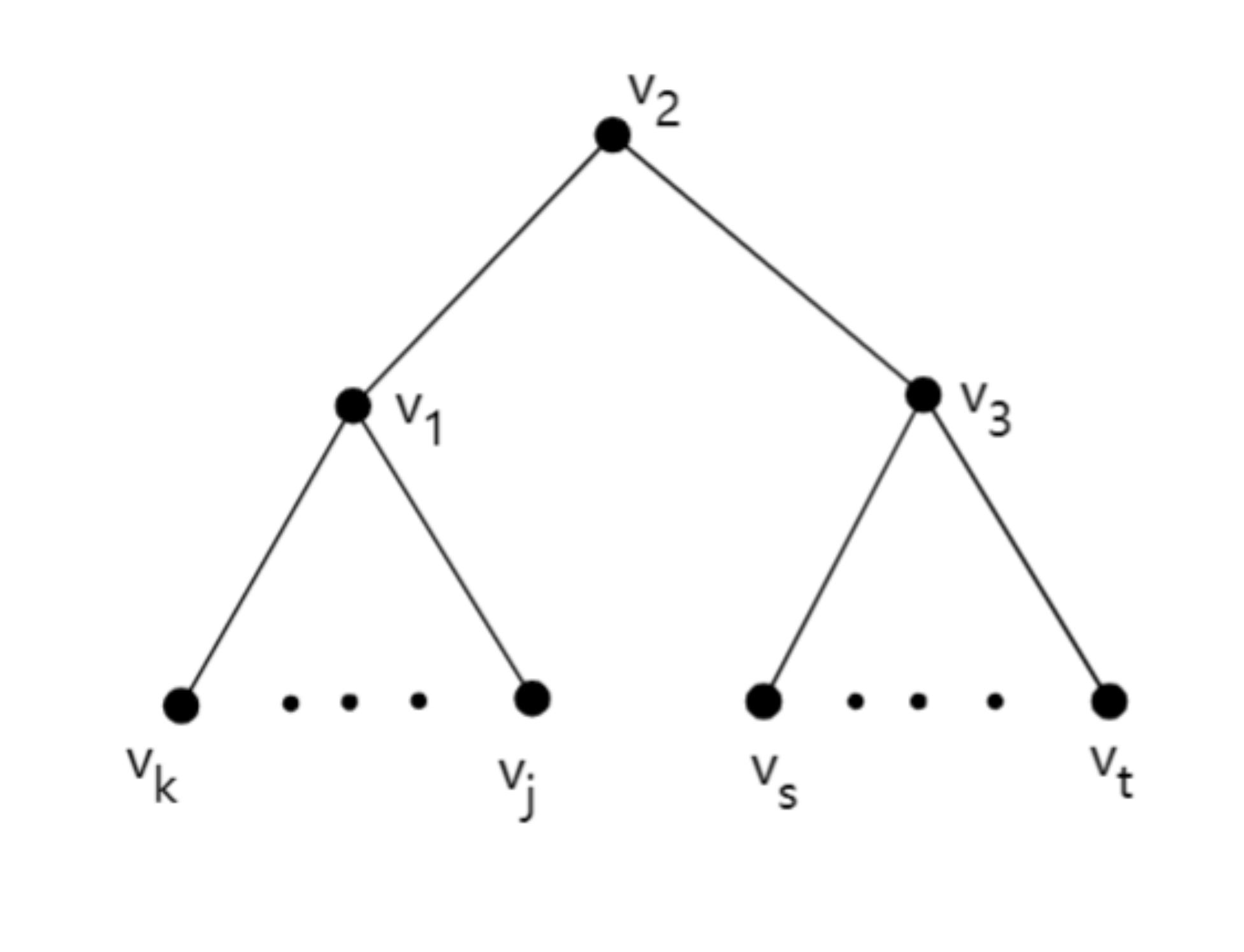}
\caption{$T(p, n-3-p) (2\le p\le n-5)$}
\label{f8}
\end{minipage}
\quad
\begin{minipage}[b]{0.45\linewidth}
\includegraphics[height=4.5cm]{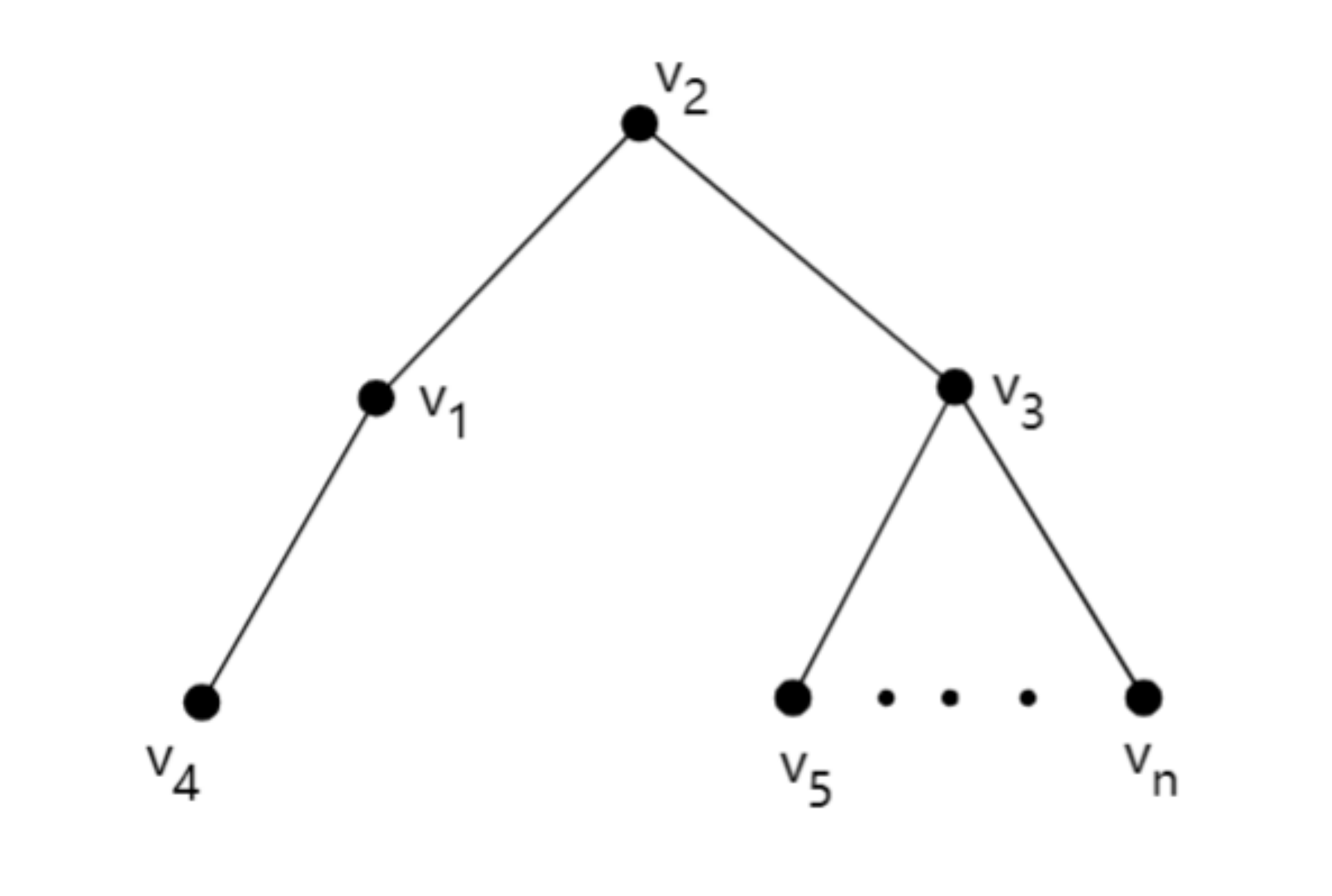}
\caption{$T(1, n-4)$}
\label{f9}
\end{minipage}
\end{figure}

If $d(v_1)\geq 3$ and $d(v_3)\geq 3$ (see Figure \ref{f8}), let $\phi(v_k), \phi(v_j)\in A_1$ and $\phi(v_s), \phi(v_t)\in B_1$, then the other leaves can be embedded in $A\cup B$ since $G$ is a balanced complete $(k+1)$-partite graph, and $T_n$ can be embedded into $G$. Now consider $d(v_1)=2$ or $d(v_3)=2$, i.e. $T_n=T(1, n-4)$. We claim that $T(1, n-4)$ cannot be embedded into $G$. Without loss of generality, assume that $d(v_1)=2$ (see Figure \ref{f9}). Suppose that $T(1, n-4)$ can be embedded into $G$ and $v_3$ is embedded into $A_1$ (without loss of generality). Then all $n-3$ vertices connected to $v_3$ should be embedded outside $A_1$. Since $|V(G)\setminus A_1|=n-3$, $\phi(v_1)$ must be in $A_1$. Since all vertices of $T(1, n-4)\setminus\{v_1, v_3\}$ are connected to $v_1$ or $v_3$, then every vertex in $T(1, n-4)\setminus\{v_1, v_3\}$ must be embedded outside $A_1$, i.e. there are $n-2$ vertices should be embedded outside $A_1$, but $|V(G)\setminus A_1|=n-3$. This is impossible. So in this case $T(1, n-4)$ cannot be embedded into $G$. This is situation (ii) in Theorem \ref{theo1}.

{\em Case 3.} $|P|=4$.

\begin{figure}[ht]
    \centering
    \includegraphics[height=4.5cm]{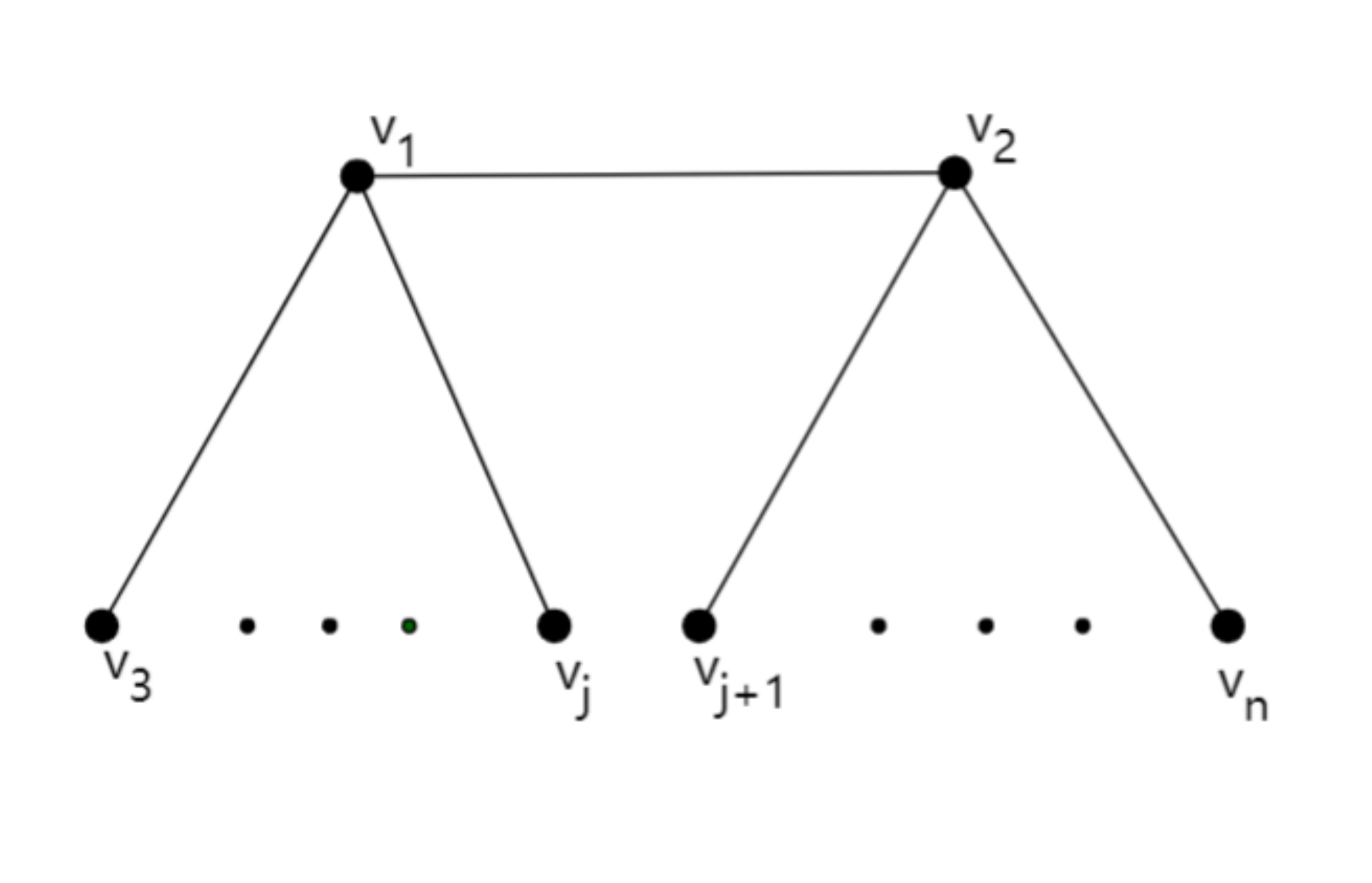}
    \caption{$T_n$ has n-2 leaves}\label{fig1}
\end{figure}

In this case we label the vertices of $T_n$ as in Figure \ref{fig1}. Note that $\{v_1, v_2, ..., v_n\}$ is a conventional labelling. Since $d(v_1)+d(v_2)=n$ and $\Delta(T_n)\leq n-3$, $d(v_1), d(v_2)\geq 3$.

{\em Case 3.1.} $\Delta(G)\geq n-1$.

Let $u, v\in V(G)$ such that $d(u)=\Delta(G)$ and $uv\in E(G)$. Embed $\phi(v_2)=u$ and $\phi(v_1)=v$. By Lemma \ref{lemma23} (iv), $\phi$ can be extended to $\{v_1, ..., v_{n-2}\}$. Recall that $v_{n-1}$ and $v_n$ have the same parent (see Figure \ref{fig1}). Since $d(u)\geq n-1$, there are two free vertices in $N(u)\setminus\{\phi(v_1), ..., \phi(v_{n-2})\}$ to embed $v_{n-1}$ and $v_n$. So $T_n$ can be embedded in $G$.

{\em Case 3.2.} $\Delta(G)=n-2$.

Let $\phi(v_2)=u_2'$ such that $d(u_2')=n-2$. We claim that there exists a vertex $u_1'\in N(u_2')$ such that there exists a vertex $u_3'\in N(u_1')\setminus N[u_2']$. Otherwise $\cup_{u\in N[u_2']}=N[u_2']$. So $N[u_2']$ is a component with $n-1$ vertices, contradicting that $G$ is a connected graph on at least $n$ vertices. Let $\phi(v_1)=u_1'$ and $\phi(v_3)=u_3'$. By Lemma \ref{lemma23} (iv), $\phi$ can be extended to $\{v_1, ..., v_{n-2}\}$. Since $d(u_2')=n-2$ and $u_3'\notin N(u_2')$, there are two free vertices in $N(u_2')\setminus\{\phi(v_1), ..., \phi(v_{n-2})\}$ to embed $v_{n-1}$ and $v_n$. So $T_n$ can be embedded into $G$.

{\em Case 3.3.} $G$ is $(n-3)$-regular.

\begin{claim} \label{case2}
If there exists a vertex $u_2'\in V(G)$ such that there exists a vertex $u_1'\in N(u_2')$ with two neighbors $u_x, u_y\in N(u_1')\setminus N[u_2']$, then $T_n$ can be embedded into $G$.
\end{claim}
{\em Proof.} Let $\phi(v_1)=u_1'$, $\phi(v_2)=u_2'$, $\phi(v_3)=u_x$ and $\phi(v_4)=u_y$. By Lemma \ref{lemma23} (iv), $\phi$ can be extended to $\{v_1, ..., v_{n-2}\}$. Since $u_2', u_x, u_y\notin N(u_2')$ and $d(u_2')=n-3$, there are two free vertices in $N(u_2')\setminus\{\phi(v_1), ..., \phi(v_{n-2})\}$ to embed $v_{n-1}$ and $v_n$. So $T_n$ can be embedded into $G$.
\q

\begin{claim} \label{case22}
If $T_n$ cannot be embedded into $G$, then for any vertex $u_x\in V(G)$, there exists a vertex $u_y\in N(u_x)$ such that $u_y$ has exactly one neighbor outside $N[u_x]$.
\end{claim}
{\em Proof.} By Claim \ref{case2}, any vertex in $N(u_x)$ has at most one neighbor outside $N[u_x]$. If $\cup_{u\in N[u_x]}N[u]\subseteq N[u_x]$, then $N[u_x]$ is a component on $n-2$ vertices, a contradiction.
\q

If $T_n$ cannot be embedded into $G$, let $u_x\in V(G)$ and $u_y\in N(u_x)$ such that $u_z\in N(u_y)\setminus N[u_x]$. Since $G$ is (n-3)-regular and $u_z\in N(u_y)\setminus N[u_x]$, there is exactly one vertex $u_a\in N(u_x)\setminus N(u_y)$ and $N(u_y)=N[u_x]\cup\{u_z\}\setminus\{u_y, u_a\}$. By Claim \ref{case2}, $N[u_y]\setminus\{u_x\}\subseteq N[u_z]$. Since $u_yu_a\notin E(G)$, $|d(u_a)|=n-3$ and $|N[u_x]|=n-2$, there is exactly one vertex $u_b\in N(u_a)\setminus N(u_x)$. If $u_b=u_z$, then $N(u_z)=N(u_x)$. For any vertex $u\in N(u_x)$, $u_z\in N(u)\setminus N(u_x)$. Then by Claim \ref{case2}, $N[u_x]\cup\{u_z\}$ is a component on $n-1$ vertices, a contradiction. So $u_b\not=u_z$. Since $u_a$ has exactly one neighbor outside $N(u_x)$, $N[u_a]=N[u_x]\cup\{u_b\}\setminus\{u_y\}$. Since all vertices in $N(u_x)$ have at most one neighbor outside $N(u_x)$ and $u_z$ is a neighbor of all vertices in $N(u_x)\setminus\{u_a\}=N(u_a)\setminus\{u_b, u_x\}\cup\{u_y\}$, by Claim \ref{case2}, $u_b$ is not adjacent to any of these vertices, i.e, $u_b$ has no neighbor in $N(u_a)\setminus\{u_b\}$, contradicting Claim \ref{case2} again. So $T_n$ can be embedded into $G$ in this case.
\q

\section{Tree-star Ramsey number }\label{sec4}
In this section, we apply Theorem \ref{theo1} to obtain the Ramsey number of a tree versus a star for some cases.

{\em Proof of Theorem \ref{cor1}.} For any red-blue-coloring of $E(K_{m+n-3})$, let $G_R$ and $G_B$ be the red graph and the blue graph respectively. If $\Delta(G_B)\geq m$, then there is a blue $K_{1, m}$. So we may assume that $\Delta(G_B)\leq m-1$. Hence $\delta(G_R)\geq (m+n-4)-(m-1)=n-3$. Since $m+n-3$ is not a linear combination of $n-1$ and $n-2$, then $G_R$ contains a component with order at least $n$. By Theorem \ref{theo1}, $T_n$ can be embedded into $G_R$.
\q

\begin{fact}\label{fact1}
For positive integers $m$ and $n$, $m+n-4$ is a linear combination of $n-1$ and $n-2$ and $m+n-3$ is not a linear combination of $n-1$ and $n-2$ if and only if $m=k(n-1)+3$ and $0\leq k\leq n-5$.
\end{fact}
{\em Proof of Fact \ref{fact1}.} Necessity: If $m+n-4=k(n-1)+l(n-2)$, $k\geq 0$ and $l\geq 1$, then $$m+n-3=k(n-1)+l(n-2)+(n-1)-(n-2)=(k+1)(n-1)+(l-1)(n-2),$$ contradicting that $m+n-3$ is not a linear combination of $n-1$ and $n-2$. Therefore $l=0$ and $m=(k-1)(n-1)+3$, where $k\geq 1$. Let us rewrite $m=k(n-1)+3$, where $k\geq 0$. If $k\geq n-4$, then $$m+n-3=k(n-1)+n=(k-n+4)(n-1)+(n-2)(n-2),$$ and $m+n-3$ is a linear combination of n-1 and n-2. So $0\leq k \leq n-5$.

Sufficiency: If $m=k(n-1)+3$ and $0\leq k\leq n-5$, then $m+n-4=(k+1)(n-1)$ is a linear combination of $n-1$ and $n-2$. If $m+n-3$ is a linear combination of $n-1$ and $n-2$, then $m+n-3=k(n-1)+n=k'(n-1)+l'(n-2)$ ($k', l'\geq 0$). So $(k+2-k')(n-1)=(l'+1)(n-2)$. Since $n-1$ is relatively prime to $n-2$, $k+2-k'=c(n-2)$ for $c\geq 1$. Therefore $k\geq (n-2)-2+k'\geq n-4$, a contradiction.
\q

For a positive integer $a$ dividing $n-3$, let $B(n-3+a)$ be the balanced complete $\frac{n-3+a}{a}$-partite graph on $n-3+a$ vertices.

{\em Proof of Theorem \ref{cor2}.} Since $m+n-4=(k+1)(n-1)$, we can partition $V(K_{m+n-4})$ into $k+1$ disjoint equal parts with $n-1$ vertices in each part. Color edges in each part red and other edges blue. Then $\Delta(G_B)\leq m+n-4-(n-1)=m-3$. So this coloring yields neither a red $T_n$ nor a blue $K_{1, m}$. Therefore $R(T_n, K_{1, m})\geq m+n-3$.

(i). By Fact \ref{fact1} and Theorem \ref{cor1}, $R(T_n, K_{1, m})\leq m+n-3$. Therefore $R(T_n, K_{1, m})=m+n-3$.

(ii). For any red-blue-coloring of $E(K_{m+n-3})$, let $G_R$ and $G_B$ be the red graph and the blue graph respectively. If $\Delta(G_B)\geq m$, then there is a blue $K_{1, m}$. So we may assume that $\Delta(G_B)\leq m-1$. Hence $\delta(G_R)\geq (m+n-4)-(m-1)=n-3$. By Fact \ref{fact1}, $G_R$ contains a component with at least $n$ vertices. By Theorem \ref{theo1}, if $m+n-3$ is not a linear combination of $n-1$, $n-2$, $n-3+a_1$, ..., $n-3+a_d$, then $T(1, n-4)$ can be embedded into $G_R$. So $R(T(1, n-4), K_{1, m})\leq m+n-3$. Therefore $R(T(1, n-4), K_{1, m})=m+n-3$.

If $m+n-3=k_1(n-1)+l_1(n-2)+s_1(n-3+a_1)+...+s_d(n-3+a_d)$, for red-blue-coloring of $E(K_{m+n-3})$, let $G_R$ and $G_B$ be the red graph and the blue graph respectively. Let $G_R=k_1K_{n-1}\cup l_1K_{n-2}\cup s_1B(n-3+a_1)\cup ...\cup s_dB(n-3+a_d)$ and $G_B=G_R^c$. By Theorem \ref{theo1}, $T(1, n-4)$ cannot be embedded into $G_R$. Since $\delta(G_R)\geq n-3$, then $\Delta(G_B)\leq (m+n-4)-(n-3)=m-1$. So there is neither a red $T_n$ nor a blue $K_{1, m}$. Therefore $R(T(1, n-4), K_{1, m})\geq m+n-2$. The same lower bound was given in \cite{GV}. So $R(T(1, n-4), K_{1, m})= m+n-2$.

(iii). For any red-blue-coloring of $E(K_{m+n-3})$, let $G_R$ and $G_B$ be the red graph and the blue graph respectively. If $\Delta(G_B)\geq m$, then there is a blue $K_{1, m}$. So we may assume that $\Delta(G_B)\leq m-1$. Hence $\delta(G_R)\geq n-3$. By Fact \ref{fact1} and $m+n-3$ is not a linear combination of $2n-6$, $n-1$ and $n-2$, there exists a component $G'$ such that $n\leq |V(G')|\not=2n-5$. By Theorem \ref{theo1}, $T(p, n-3-p)$ can be embedded into $G'$. Therefore $R(T(p, n-3-p), K_{1, m})= m+n-3$. If $m+n-3=r(2n-6)+s(n-1)+t(n-2)$, then we can partition $K_{m+n-3}$ into $r+s+t$ disjoint parts such that $r$ parts have order $2n-6$ and $s$ parts have order $n-1$ and $t$ parts have order $n-2$. Let $G_R=rK_{n-3, n-3}\cup sK_{n-1}\cup tK_{n-2}$ and $G_B=G_R^c$. Then $\Delta(G_B)\leq m+n-4-(n-3)=m-1$ and $T(p, n-3-p)$ cannot be embedded into $G$. Therefore $R(T(p, n-3-p), K_{1, m})\geq m+n-2$. The same lower bound was given in \cite{GV}. So $R(T(1, n-4), K_{1, m})= m+n-2$.
\q

\bigskip

{\bf Acknowledgement}
The research is supported in part by National Natural Science Foundation of China (No. 11931002).

\end{document}